\documentclass[a4 paper,12pt,bezier]{article}
\usepackage[top=3cm,bottom=3cm,left=2.5cm,right=2.5cm]{geometry}
\usepackage{amsmath}
\usepackage{amsfonts,amsthm,amssymb}
\usepackage{amsfonts}
\usepackage{graphics,subfigure,caption2}
\usepackage{tikz, color}
\usepackage{graphics,epstopdf}
\usepackage{latexsym,bm}
\usepackage{amsfonts,amsthm,amssymb,mathrsfs,bbding}
\usepackage{CJK}
\usepackage{color}
\usepackage{indentfirst}
\newtheorem{lem}{Lemma}[section]
\newtheorem{thm}[lem]{Theorem}
\newtheorem{cor}[lem]{Corollary}

\newtheorem{cl}{Claim}
\theoremstyle{plain}

\usepackage{cite}

\newcounter{countclaim}

\def \clproof{\noindent {\it Proof}. }

\newcommand {\clproofend}
{{\hfill$\natural$}}
\begin{document}

\begin{CJK}{GBK}{song}
\title{\text{Chv\'{a}tal-Erd\H{o}s} condition for 2-factors with at most two components in graphs}
  \author{Tao Tian$^1$, Liming Xiong$^{2,}$\footnote{Corresponding author. E-mail: taotian0118@163.com (T. Tian), lmxiong@bit.edu.cn (L. Xiong), weigenyan@263.net (W. Yan)}, Weigen Yan$^{3}$\\
\small $^1$  School of Mathematics and Statistics, Key Laboratory of Analytical Mathematics and\\  \small  Applications $($Ministry of Education$)$, Fujian Key Laboratory of Analytical Mathematics and\\ \small  Applications (FJKLAMA), Center for Applied Mathematics of Fujian Province $($FJNU$)$,\\ \small Fujian Normal University,
  \small Fuzhou 350117, China.\\
 \small $^2$ School of Mathematics and Statistics, Beijing Institute of Technology, Beijing 100081,  China\\
 \small $^3$  School of Science, Jimei University, Xiamen 361021, China}

\date{}
\maketitle
\maketitle {\flushleft\bf Abstract}: It is well-known that Chv\'{a}tal and Erd\H{o}s stated that any graph of order at least three whose independence number is no greater than its connectivity is Hamiltonian; that any graph whose independence number is no greater than its connectivity minus one is Hamilton-connected; and that any graph whose independence number is no greater than its connectivity plus one is traceable. Kaneko and Yoshimoto [J. Graph Theory 43 (2003) 269--279] showed that every 4-connected graph of order at least six has a 2-factor with two components if its independence number is no greater than its connectivity. In this paper, we show that any connected graph of order at least three times its connectivity plus three has a 2-factor with at most two components, except for one exceptional class, if its independence number is no greater than its connectivity plus one. Our result is best possible.
\maketitle {\flushleft\textit{\bf Keywords}: Chv\'{a}tal-Erd\H{o}s condition, 2-factor, number of components, graph}

\section{ Introduction }

We follow Bondy and Murty \cite{1} for undefined terms and notation. We consider finite, undirected and loopless graphs only. Let $G$ be a graph with vertex set $V(G)$ and edge set $E(G$.  For a vertex $x$ of $G$, denote by $N_{G}(x)$  the neighborhood of $x$ in $G$, and by $d_{G}(x)$  the degree of $x$ in $G$. Let $S,T$ be two disjoint vertex sets and $S,T\subseteq V(G)$. We define $N_{G}(S)=\cup_{x \in S}N_{G}(x)$ and denote by $E[S,T]$ the set of edges of $G$ with one end in $S$ and the other end
in $T$.  Let $G_1$ and $G_2$ be two disjoint graphs. Let $G_1\cup G_2$  and $G_1 \vee G_2$ denote the union and the join of $G_1$ and $G_2$, respectively.  The notation $kG$ is short for $\underset{k}{\underbrace{G\cup G\cup...\cup G}}$.

Let $H\subseteq G$ denote that  $H$ is a subgraph of $G$. We denote the number of components of $H$ by $\omega(H)$. A {\sl spanning subgraph} of a graph $G$ is a subgraph with vertex set  $V(G)$ and with edge set  obtained from $E(G)$ by deleting some edges. A graph is called {\sl Hamiltonian} (resp. {\sl traceable}) if it contains a spanning cycle (resp. {\sl spanning path}). A graph is called {\sl Hamilton-connected} if each pair of vertices are connected by a spanning path. A subgraph of $G$ is called a {\sl 2-factor} if it is a spanning subgraph of $G$ with each vertex of degree 2.

Let $C_{1},C_{2},\ldots,C_{k}$ be cycles in a graph $G$. We say that $G$ is covered by $C_{1},C_{2},\ldots,C_{k}$  if
 $\cup _{i=1}^k V(C_{i})=V(G)$. Furthermore, if $C_{1},C_{2},\ldots,C_{k}$  are disjoint from each other, then $C_{1},C_{2},\ldots,C_{k}$  form a 2-factor of $G$. Obviously, the property that
$G$ has a 2-factor with at most two components is stronger than having a spanning path and weaker than having a spanning cycle, i.e.,

Hamiltonicity $\Rightarrow$  2-factors with at most two
components $\Rightarrow$ Traceability.

In terms of computational complexity, the problem of deciding whether a given graph is Hamiltonian (or traceable) is generally NP-complete. This implies that deciding whether a given graph has a 2-factor with at most two components is also an NP-complete problem.

Denote by $\alpha(G)$ and $\kappa (G)$ the independence number and connectivity of $G$, respectively. In 1972, Chv\'{a}tal and Erd\H{o}s gave  sufficient conditions for a graph to be Hamiltonian, traceable, and Hamilton-connected, respectively, as follows.

\begin{thm}{\rm (Chv\'{a}tal and Erd\H{o}s \cite{3})}\label {Th2}
Let $G$ be a connected graph of order at least $3$. Then each of the following holds.

\item(a) If $\alpha(G)\leq \kappa(G)$, then $G$  is Hamiltonian.
\item(b) If $\alpha(G)\leq \kappa(G)+1$, then $G$  is traceable.
\item(c) If $\alpha(G)\leq \kappa(G)-1$, then $G$  is Hamilton-connected.

\end{thm}

Obviously, a graph $G$ is Hamiltonian if and only if a graph $G$ can be covered by exactly one cycle. In \cite{24},  Kouider proved the following beautiful result related  the  Chv\'{a}tal and Erd\H{o}s condition, which was first posed by Fournier \cite{31}.

\begin{thm} {\rm (Kouider  \cite{24})} \label {Th4}
Every 2-connected graph $G$  can be covered by $\lceil \alpha(G)/ \kappa(G) \rceil$ cycles.

\end{thm}

Particularly, when $\alpha(G)= k+1$ ($k\geq2$), Amar et al. obtained the following  result, which  is itself useful.

\begin{thm}{\rm (Amar et al. \cite{2})} \label{Th6}  Let $G$ be a $k$-connected simple graph with $\alpha(G)= k+1$ ($k\geq2$). If $C$ is a longest cycle of $G$, then  $G-C$ is a complete graph.
\end{thm}

Note that even if a graph  $G$ can be covered by  some cycles, these cycles may not form a 2-factor of $G$. In \cite{5}, Kaneko and Yoshimoto proved a sufficient Chv\'{a}tal and Erd\H{o}s condition implying that a graph has a 2-factor with small number of components, as follows.

\begin{thm}{\rm (Kaneko and Yoshimoto  \cite{5})}\label {Th3}
Let $G$ be a 4-connected graph with $|V(G)|\geq 6$. If $\alpha(G)\leq \kappa(G)$, then $G$ has a 2-factor with two components.

\end{thm}


It is known that degree conditions are classic conditions for a graph admitting a Hamiltonian cycle. Similarly, sufficient Dirac-type or Ore-type degree conditions for graphs containing a 2-factor with exactly $k$-components have yielded many significant research results; see, for example, \cite{6,7,8,9,13,32}. In addition,  sufficient conditions for graphs containing 2-factor with  bounded number of components also attracted scholars' attention; see, for example, \cite{10,11,15,25,17,33}. Other works related to Chv\'{a}tal--Erd\H{o}s condition can be found in \cite{23,16,29,12,26,27}. For other works related to 2-factors, we refer the reader to  \cite{14,19,20,21,28,22}.

Motivated by the  above results, we obtain the main result of this paper as follows. Before stating it, we first define  a  class of graphs.

Let $K_{a}$ be a clique of order $a$, i.e., a complete graph with  $a$ vertices. Define  $$\mathcal {G}=\{G|G_{1}\subseteq G \subseteq G_{2}\}\cup \{G_{3}\},$$ where $G_{1}$ is obtained from two cliques $K_{n-2}$ and $K_{2}$ by connecting two vertices of $K_{n-2}$ to a same vertex in $K_{2}$; $G_{2}$ is obtained from two cliques $K_{n-1}$ and $K_{1}$ by connecting one vertex of $K_{n-1}$ to the vertex of $K_{1}$; and $G_{3}$ is obtained from two cliques $K_{n-2}$ and $K_{2}$ by connecting one vertex of $K_{n-2}$ to exactly one vertex in  $K_{2}$, as depicted in Fig. 1.

\begin{figure}[htbp]
  \centering
  \includegraphics[width=0.8\linewidth]{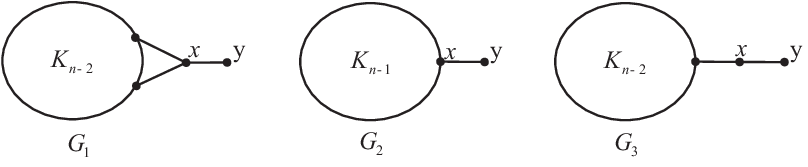}
  \caption{The graphs $G_1$,$G_2$ and $G_3$.}
\end{figure}

\begin{thm}\label {Th1}
Let $G$ be a connected graph of order at least $3\kappa(G)+3$. If $\alpha(G)\leq \kappa(G)+1$, then $G$ has a 2-factor with at most two components unless $G\in \mathcal {G}$.
\end{thm}

 \par The following three corollaries can be derived from Theorem \ref{Th1} immediately.
\begin{cor}\label {C1}
Let $G$ be a connected graph of order at least $3\kappa(G)+3$. If $\alpha(G)\leq \kappa(G)+1$, then $G$ is traceable.
\end{cor}

\begin{cor}\label {C2}
Let $G$ be a connected graph with order at least $3\kappa(G)+3$ and minimum degree at least 2. If $\alpha(G)\leq \kappa(G)+1$, then $G$ has a 2-factor with at most two components.
\end{cor}

\begin{cor}\label {C3}
Let $G$ be a connected non-Hamiltonian graph with order at least $3\kappa(G)+3$ and minimum degree at least 2. If $\alpha(G)\leq \kappa(G)+1$, then $G$ has a 2-factor with exactly two components.
\end{cor}






\section{Proof of Theorem \ref{Th1}}

\noindent The proof is by contradiction.  Suppose  $G$ has no 2-factor with at most two components. If $\alpha(G)\leq \kappa(G)$, then by Theorem \ref{Th2}, $G$  is Hamiltonian, a contradiction.
Then

\begin{equation}
	\alpha(G)=\kappa(G)+1.
\end{equation}

If $\kappa(G)=1$, then $G$ has a cycle; otherwise, by $\alpha(G)=2$, $G$ is connected acyclic graph with at most two vertices of degree 1, i.e., $G$ is a spanning path with at least six vertices, we can easily find an independent set of $G$ with  cardinality at least 3, a contradiction. In what follows, let $G$ be a graph of order $n$, and let $C$ be a longest cycle in $G$.

\begin{cl} \label {cl14}
 If $\kappa(G)=1$, then $G-C$ is a complete graph.

\end{cl}
\clproof  We set $V(C)=\{v_{1},v_{2},\ldots, v_{l}\}$, where $C=v_{1}v_{2}\cdots v_{l-1}v_{l}v_{1}$ and $l\geq3$. By $\alpha(G)=2$, $\omega(G-C)\leq2$.   Without loss of generality, we assume that $x\in V(G-C)$ and $xv_{1}\in E(G)$. Then $xv_{2},xv_{l}\notin E(G)$; otherwise,  $C'=xv_{2}\cdots v_{l-1}v_{l}v_{1}x$ or $C''=xv_{1}v_{2}\cdots v_{l-1}v_{l}x$ is a cycle of length $|V(C)|+1$, a contradiction. Suppose first that $G-C$ has two components, $H_{1}$ and $H_{2}$.  Without loss of generality, we further assume that $x\in V(H_{1})$  and $y\in V(H_{2})$. By $\alpha(G)=2$, $v_{2}y,v_{l}y\in E(G)$. Similarly, we have $yv_{3}\notin E(G)$. Then, by $\alpha(G)=2$, $xv_{3}\in E(G)$. However, $v_{1}xv_{3}\cdots v_{l}yv_{2}v_{1}$
is a cycle of length $|V(C)|+2$, a contradiction.

 Then $G-C$ has exactly one component $H_{1}$. Suppose that $H_{1}$ is not complete. Then, there exists one pair of vertices $y,z\in V(H_{1})$ and $yz\notin E(H_{1})$. By $\alpha(G)=2$ and $xv_{l}\notin E(G)$, we can choose one vertex $u\in \{y,z\}\setminus \{x\}$ such that $uv_{l}\in E(G)$. Then $v_{1}v_{2}\cdots v_{l}uPxv_{1}$ is a cycle of length at least $|V(C)|+2$, where $P$ is a path in $H_{1}$ connecting vertices $u$ and $x$, a contradiction.
Thus,  Claim \ref{cl14} holds. \clproofend

By Claim \ref{cl14} and Theorem \ref{Th6}, $G-C$ is a complete graph. Denote by $H$  the component of $G-C$. Then $1\leq|V(H)|\leq2$. Otherwise, $G$ has a 2-factor $C\cup C_{H}$, where $C_{H}$ is a spanning cycle of $H$, a contradiction. Then $H=K_{1}$ or $K_{2}$. Throughout the proof, we set $V(H)=\{x\}$ if $H=K_{1}$,  and  $V(H)=\{x,y\}$ if $H=K_{2}$.  In the following, we distinguish two cases $\kappa(G)=1$ and $\kappa(G)\geq 2$.

\par {\bfseries Case 1:}\label{cas1}\ $\kappa(G)=1$.

Suppose first that $H=K_{1}$. Let $V(C)=\{v_{1},v_{2},\ldots, v_{n-1}\}$, where $C=v_{1}v_{2}\cdots v_{n-1}v_{1}$. Since $\kappa(G)=1$, $|E(V(C),\{x\})|=1$. Without loss of generality, we assume that $xv_{n-1}\in E(G)$ and $xv_{i}\notin E(G)$ for $i=1,2,\ldots,n-2$. By $\alpha(G)=2$, $G[\{v_{1},v_{2},\ldots, v_{n-2}\}]\cong K_{n-2}$ and $2\leq |E[\{v_{1},v_{2},\ldots, v_{n-2}\},\{v_{n-1}\}]|\leq n-2$. Then $G_{1}\subseteq G\subseteq G_{2}$, where $G_{1}$ and $G_{2}$ are depicted in Fig. 1.


Now suppose  that $H=K_{2}$. Let $V(C)=\{v_{1},v_{2},\ldots, v_{n-2}\}$, where $C=v_{1}v_{2}\cdots v_{n-2}v_{1}$. Since $\kappa(G)=1$, without loss of generality, we assume that $xv_{n-2}\in E(G)$ and $yv_{i}\notin E(G)$ for $i=1,2,\ldots,n-3$.
If $yv_{n-2}\in E(G)$, then $xv_{i}\notin E(G)$ for $i=1,2,\ldots,n-3$ by $\kappa(G)=1$. Then $G[\{v_{1},v_{2},\ldots, v_{n-3}\}]\cong K_{n-3}$ by $\alpha(G)=2$. However,  $v_{1}v_{2}\cdots v_{n-3}v_{1}$ and $v_{n-2}xyv_{n-2}$ constitute a 2-factor of $G$ and the  number of components of the 2-factor is exactly 2 at this case, a contradiction. Then $yv_{n-2}\notin E(G)$ and so $E[V(C),\{y\}]=\emptyset$. By $\alpha(G)=2$, $G[\{v_{1},v_{2},\ldots, v_{n-2}\}]\cong K_{n-2}$; otherwise, one pair of nonadjacent vertices in $\{v_{1},v_{2},\ldots, v_{n-2}\}$ and the vertex $y$ constitute an independent set $G$ with cardinality 3, a contradiction. Then $xv_{i}\notin E(G)$ for $i=1,2,\ldots,n-3$; otherwise, we can easily find  a cycle $C'$ containing vertex $x$ with $V(C)\subset V(C')$, a contradiction. Then, the graph $G\cong G_{3}$, where $G_{3}$ is depicted in Fig. 1.


\par {\bfseries Case 2:}\label{cas2}\ $\kappa(G)\geq 2$.

We assume that the  cycle $C$  has a fixed orientation. Let $u$ and $v$ be two vertices in $C$. We denote by $u^{+}$ and $u^{-}$, respectively,  the
successor and the predecessor of $u$ in $C$, by $u^{++}$ the successor of $u^{+}$, by $u^{t+}$ the $t$-th successor of $u$, by $u^{--}$ the predecessor of $u^{-}$,
by $u\overrightarrow{C}v$ (resp. $u\overleftarrow{C}v$) the subpath in $C$ from $u$ to $v$  that contains $u+$ (resp. $u^{-}$),
and by $\overrightarrow{C}[u,v]$ the section of $C$ that contains $u$ and $v$.  Particularly, $u\overrightarrow{C}v=u$ if $u=v$.

Denote by $u_{1},u_{2},\ldots,u_{k}$ the neighborhoods of $H$ on $C$, i.e., $N_{C}(H)=N_{G}(V(H))\cap V(C)=\{u_{1},u_{2},\ldots,u_{k}\}$, where the subscripts are labeled in order according to the orientation of the cycle $\overrightarrow{C}$. Denote by $U=\{u_{1},u_{2},\ldots,u_{k}\}$, $U^{+}=\{u_{1}^{+},u_{2}^{+},\ldots,u_{k}^{+}\}$ and
 $U^{-}=\{u_{1}^{-},u_{2}^{-},\ldots,u_{k}^{-}\}$.
 Since $G$ is $\kappa(G)$-connected, $k\geq \kappa(G)$. In the following, the subscripts are taken modulo $k$.

 The following Claims \ref{cl1}
 and \ref{cl2} have been used in the proof of Chv\'{a}tal-Erd\H{o}s Theorem and will be repeatedly used in our later proof.

  \begin{cl}\rm{([2])} \label {cl1}
 For any $u_{i}\in U$, $u_{i}u_{i+1}\notin E(C)$.

\end{cl}


 \begin{cl}\rm{([2])}\label {cl2}
 For any $x\in V(H)$, $\{x\}\cup U^{+}$ and $\{x\}\cup U^{-}$ are independent sets of $G$.

\end{cl}



Since $k\geq \kappa(G)$,  by Claim \ref{cl2} and (2), we have $k=\kappa(G)$.

\begin{cl} \label {cl3}
 For any $u_{i},u_{j}\in U$ $(i\neq j)$, there exists a path $u_{i}Pu_{j}$ in $G$, where $P$ is a path in $H$ connecting $u_{i}$ and $u_{j}$ such that $V(P)=V(H)$.

\end{cl}
\clproof If $H=K_{1}$, then by $N_{C}(H)=U$, the claim holds immediately. Now, we consider the case when $H=K_{2}$. Then $xy\in E(H)$. Suppose that there does not exist a path $P$ in $H$ connecting $u_{i}$ and $u_{j}$ such that $V(P)=V(H)$. Since $u_{i},u_{j}\in N_{C}(H)$, without loss of generality, we can assume that $x\in V(H)$ satisfying $u_{i}x,u_{j}x\in E(H)$.  Then $u_{i}y,u_{j}y\notin E(H)$; otherwise, $u_{i}yxu_{j}$ or $u_{i}xyu_{j}$ is a path of $G$  containing all vertices of $H$, a contradiction. Then, by $N_{C}(H)=U$, $d_{G}(y)\leq |U|-1=\kappa(G)-1$, a contradiction. Thus,  Claim \ref{cl3} holds. \clproofend

Denote by $t_{i}=|V(\overrightarrow{C}[u_{i}^{+},u_{i+1}^{-}])|$. By Claim \ref{cl1}, $t_{i}\geq1$. For convenience, we set $T=\{t_{1},t_{2},\ldots,t_{\kappa(G)}\}$, $T_{1,2}=\{t_{i}\mid t_{i}\in T~\text{and}~1\leq t_{i}\leq2\}$, $T_{3}=\{t_{i}\mid t_{i}\in T~\text{and}~t_{i}=3\}$, and $T_{\geq3}=\{t_{i}\mid t_{i}\in T~\text{and}~t_{i}\geq 3\}$. Then $T=T_{1,2}\cup T_{\geq3}$.


\begin{cl} \label {cl4}
 $T_{\geq3}\neq \emptyset.$

\end{cl}
\clproof We establish the claim by contradiction. Suppose that $T_{\geq3}= \emptyset.$ Then $T=T_{1,2}$. Then $|V(G)|\leq |U|+2|T_{1,2}|+|V(H)|\leq 3\kappa(G)+2$, contradicting the condition that $|V(G)|\geq 3\kappa(G)+3$. \clproofend

\begin{cl} \label {cl13}
  If $t_{i}\geq3$,  then $u_{i}^{+}u_{j}^{-}\notin E(G)$ for any $j\in \{1,2,\dots,\kappa(G)\}\setminus \{i\}$.
\end{cl}
\clproof  We establish the claim by contradiction. Suppose that  $u_{i}^{+}u_{j}^{-}\in E(G)$ for some $j\in \{1,2,\dots,\kappa(G)\}\setminus \{i\}$. Then, by Claim \ref{cl3}, $u_{i}^{+}\overrightarrow{C}u_{j}^{-}u_{i}^{+}$ and $u_{i}Pu_{j}\overrightarrow{C}u_{i}$ constitute a 2-factor of $G$ and the  number of components of the 2-factor is  exactly 2, a contradiction. \clproofend

\begin{cl} \label {cl5}
 If $t_{i}\geq 3$, then $u_{i}^{-}u_{i}^{+},u_{i+1}^{-}u_{i+1}^{+}\in E(G)$.

\end{cl}
\clproof  By Claim \ref{cl13},  $u_{i}^{+}u_{j}^{-}\notin E(G)$ for any $j\in \{1,2,\dots,\kappa(G)\}\setminus \{i\}$. Then $u_{i}^{-}u_{i}^{+}\in E(G)$; otherwise, by Claims \ref{cl1} and \ref{cl2}, $\{u_{i}^{+},x\}\cup U^{-}$ is an independent set of $G$ with cardinality $\kappa(G)+2$, a contradiction. Similarly, we can prove that $u_{i+1}^{-}u_{i+1}^{+}\in E(G)$. Thus,  Claim \ref{cl5} holds. \clproofend

\begin{cl} \label {cl6}
 If $t_{i}\geq 4$, then $u_{i}^{-}u_{i}^{++},u_{i+1}^{--}u_{i+1}^{+}\in E(G)$.
\end{cl}
\clproof  By Claim \ref{cl5}, $u_{i}^{-}u_{i}^{+},u_{i+1}^{-}u_{i+1}^{+}\in E(G)$.  For any $j\in \{1,2,\dots,\kappa(G)\}\setminus \{i\}$, $u_{i}^{++}u_{j}^{-}\notin E(G)$; otherwise, by Claim \ref{cl3}, $u_{i}^{++}\overrightarrow{C}u_{j}^{-}u_{i}^{++}$ and $u_{i}Pu_{j}\overrightarrow{C}u_{i}^{-}u_{i}^{+}u_{i}$ constitute a 2-factor of $G$ and the  number of components of the 2-factor is  exactly 2, a contradiction. Then $u_{i}^{-}u_{i}^{++}\in E(G)$; otherwise, by Claims \ref{cl1} and \ref{cl2}, $\{u_{i}^{++},x\}\cup U^{-}$ is an independent set of $G$ with cardinality $\kappa(G)+2$, a contradiction. Similarly, we can prove that $u_{i+1}^{--}u_{i+1}^{+}\in E(G)$. Thus,  Claim \ref{cl6} holds. \clproofend

\begin{cl} \label {cl7}
 If $t_{i}\geq 5$, then $u_{i}^{-}u_{i}^{l+}\in E(G)$ for all $l$ with $3\leq l \leq t_{i}-2$.
\end{cl}
\clproof  We prove the claim by induction on $l$. By Claims \ref{cl5} and \ref{cl6}, $u_{i}^{-}u_{i}^{+},u_{i+1}^{-}u_{i+1}^{+},$
 $u_{i}^{-}u_{i}^{++},u_{i+1}^{--}u_{i+1}^{+}\in E(G)$.
If $l=3$, then  $u_{i}^{3+}u_{j}^{-}\notin E(G)$ for any $j\in \{1,2,\dots,\kappa(G)\}\setminus \{i\}$; otherwise, by Claim \ref{cl3}, $u_{i}^{3+}\overrightarrow{C}u_{j}^{-}u_{i}^{3+}$ and $u_{i}Pu_{j}\overrightarrow{C}u_{i}^{-}u_{i}^{++}\overleftarrow{C}u_{i}$ constitute a 2-factor of $G$ and the  number of components of the 2-factor is  exactly 2, a contradiction. Then $u_{i}^{-}u_{i}^{3+}\in E(G)$; otherwise, by Claims \ref{cl1} and \ref{cl2}, $\{u_{i}^{3+},x\}\cup U^{-}$ is an independent set of $G$ with cardinality $\kappa(G)+2$, a contradiction. In the following, we suppose that $u_{i}^{-}u_{i}^{l+}\in E(G)$ holds for all $l$ with $3\leq l <k\leq t_{i}-2$. Then  $u_{i}^{k+}u_{j}^{-}\notin E(G)$ for any $j\in \{1,2,\dots,\kappa(G)\}\setminus \{i\}$; otherwise, by Claim \ref{cl3}, $u_{i}^{k+}\overrightarrow{C}u_{j}^{-}u_{i}^{k+}$ and $u_{i}Pu_{j}\overrightarrow{C}u_{i}^{-}u_{i}^{(k-1)+}\overleftarrow{C}u_{i}$ constitute a 2-factor of $G$ and the  number of components of the 2-factor is  exactly 2, a contradiction. Then $u_{i}^{-}u_{i}^{k+}\in E(G)$; otherwise, by Claims \ref{cl1} and \ref{cl2}, $\{u_{i}^{k+},x\}\cup U^{-}$ is an independent set of $G$ with cardinality $\kappa(G)+2$, a contradiction. Thus,  Claim \ref{cl7} holds. \clproofend


\begin{cl} \label {cl8}
 $T_{\geq3}=T_{3}$.
\end{cl}
\clproof  We establish the claim by contradiction. Suppose that there exists a $t_{i}\in T_{\geq3}$ and $t_{i}\geq 4$. Then, by Claims \ref{cl3}, \ref{cl5}--\ref{cl7}, $u_{i}Pu_{i+1}u_{i+1}^{-}u_{i+1}^{--}u_{i+1}^{+}\overrightarrow{C}u_{i}^{-}u_{i}^{(t_{i}-2)+}\overleftarrow{C}u_{i}$
is a spanning cycle of $G$, a contradiction. Thus,  Claim \ref{cl8} holds. \clproofend

\begin{cl} \label {cl9}
 $H=K_{1}$.
\end{cl}
\clproof  We establish the claim by contradiction. Suppose that $H=K_{2}$. By Claims \ref{cl4} and \ref{cl8}, there exists a $t_{i}\in T_{3}$. By Claim \ref{cl5}, $u_{i}^{-}u_{i}^{+},u_{i+1}^{-}u_{i+1}^{+}\in E(G)$. Then, by Claim \ref{cl3}, $C'=u_{i}Pu_{i+1}u_{i+1}^{-}u_{i+1}^{+}\overrightarrow{C}u_{i}^{-}u_{i}^{+}u_{i}$
is a  cycle of length $|V(C)|+1$, a contradiction. Thus,  Claim \ref{cl9} holds. \clproofend

\begin{cl} \label {cl10}
 $|T_{3}|\geq2$.
\end{cl}
\clproof  We establish the claim by contradiction. Suppose that $|T_{3}|\leq1 $. Then, by Claim \ref{cl9}, $|V(G)|\leq |U|+2|T_{1,2}|+2=3\kappa(G)+2$, contradicting the condition that $|V(G)|\geq 3\kappa(G)+3$. Thus,  Claim \ref{cl10} holds. \clproofend

\begin{figure}
  \centering
  \includegraphics[width=0.95\linewidth]{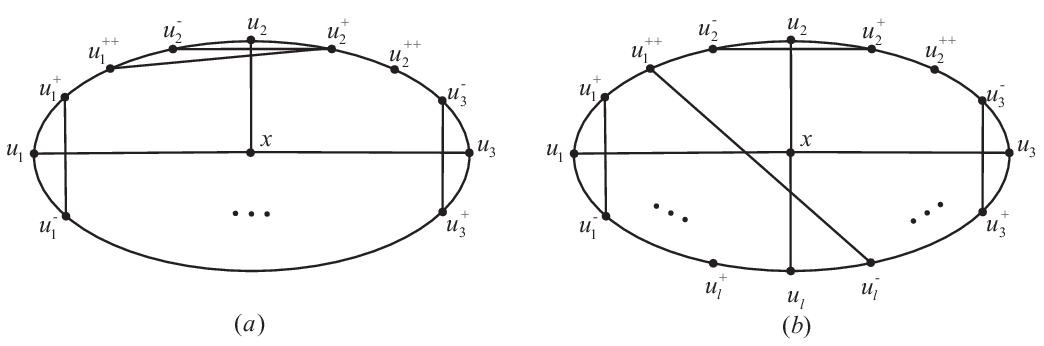}
  \caption{(a) The case when $u_{2}^{+}u_{1}^{++}\in E(G)$. (b) The case when $u_{1}^{++}u_{l}^{-}\in E(G)$.}
\end{figure}

\begin{cl} \label {cl11}
For any $t_{i},t_{j}\in T_{3}$ ($i\neq j$), $|i-j|\not\equiv 1~(\text{mod}~\kappa(G))$.
\end{cl}

\clproof We establish the claim by contradiction. Suppose that $|i-j|\equiv1~(\text{mod}~\kappa(G))$. Without loss of generality, we assume that $i=1$ and $j=2$. Then, by Claim \ref{cl5}, $u_{1}^{-}u_{1}^{+},u_{2}^{-}u_{2}^{+}\in E(G)$. By Claim \ref{cl13}, $u_{2}^{+}u_{l}^{-}\notin E(G)$ for any $l\in \{1,2,\dots,\kappa(G)\}\setminus \{2\}$. Then $u_{2}^{+}u_{1}^{++}\notin E(G)$;  otherwise,  $u_{1}xu_{2}u_{2}^{-}u_{1}^{++}u_{2}^{+}\overrightarrow{C}u_{1}^{-}u_{1}^{+}u_{1} $ is a spanning cycle of $G$ (see Fig. 2(a)), a contradiction.
Then $u_{1}^{++}u_{l}^{-}\notin E(G)$ for any $l\in \{1,2,\dots,\kappa(G)\}\setminus \{2\}$; otherwise,  $u_{1}u_{1}^{+}u_{1}^{++}u_{l}^{-}\overleftarrow{C}u_{2}^{+}u_{2}^{-}u_{2}xu_{l}\overrightarrow{C}u_{1}$ is a spanning cycle of $G$ (see Fig. 2(b)), a contradiction. Then, by Claims \ref{cl1} and \ref{cl2}, $\{u_{2}^{+},u_{1}^{++},x\}\cup (U^{-}\setminus\{u_{2}^{-}\})$ is an independent set of $G$ with cardinality $\kappa(G)+2$, a contradiction. Thus,  Claim \ref{cl11} holds. \clproofend


\begin{figure}[htbp]
  \centering
  \includegraphics[width=0.95\linewidth]{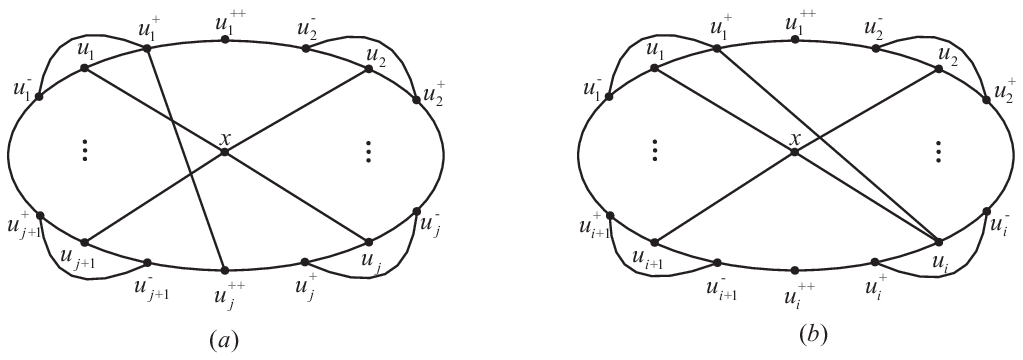}
  \caption{(a) The case when $u_{1}^{+}u_{j}^{++}\in E(G)$. (b) The case when $j=i$.}
\end{figure}

By Claim \ref{cl10}, $|T_{3}|\geq2$. Without loss of generality, we assume that $t_{1},t_{i}\in T_{3}$. By Claim  \ref{cl11}, $\kappa(G)\geq4$ and $i\in \{3,4,\ldots,\kappa(G)-1\}$. By Claims \ref{cl5}, $u_{1}^{-}u_{1}^{+},u_{2}^{-}u_{2}^{+},u_{i}^{-}u_{i}^{+},$
$u_{i+1}^{-}u_{i+1}^{+}\in E(G)$. Denote by $U^{++}=\{u_{j}^{++}\mid t_{j}\in T_{3}~\text{and}~j\in \{3,4,\ldots,\kappa(G)-1\}\}$. Then $u_{1}^{+}u_{j}^{++}\notin E(G)$, where $u_{j}^{++}\in U^{++}$; otherwise, by Claim \ref{cl5}, $u_{1}^{+}\overrightarrow{C}u_{j}^{++}u_{1}^{+}$ and $u_{1}xu_{j+1}u_{j+1}^{-}u_{j+1}^{+}\overrightarrow{C}u_{1}$ constitute a 2-factor of $G$ and the  number of components of the 2-factor is  exactly 2 (see Fig. 3(a)), a contradiction. By Claim \ref{cl13}, $u_{1}^{+}u_{j}^{-}\notin E(G)$ for $j\in \{2,\dots,\kappa(G)\}$. Then $u_{1}^{+}u_{j}\notin E(G)$, where $ j\in\{2,i,i+1\}$; otherwise, by Claim \ref{cl5}, $u_{1}^{+}\overrightarrow{C}u_{j}^{-}u_{j}^{+}\overrightarrow{C}u_{1}xu_{j}u_{1}^{+} $ is a spanning cycle of $G$ (see Fig. 3(b) for an example when $j=i$), a contradiction. By Claim \ref{cl2}, $U^{+}$ is an independent set of $G$. Then
\begin{equation}
  \begin{aligned}
d_{G}(u_{1}^{+})\leq& |N_{G}(u_{1}^+)\cap U|+|N_{G}(u_{1}^+)\cap U^-|+|N_{G}(u_{1}^+)\cap U^{++}|+|N_{G}(u_{1}^+)\cap V(H)|\\
\leq&|U\setminus\{u_{2},u_{i},u_{i+1}\}|+|\{u_{1}^{-}\}|+|\{u_{1}^{++}\}|+0 \\
=&\kappa(G)-1.
  \end{aligned}
 \end{equation}
 
However, since $G$ is $\kappa(G)$-connected, $d_{G}(u_{1}^{+})\geq \kappa(G)$, contradicting (2).

This completes the proof of Theorem \ref{Th1}.  \hfill$\blacksquare$ 

\section{Sharpness of Theorem \ref{Th1}}

\begin{figure}[htbp]
  \centering
  \includegraphics[width=0.55\linewidth]{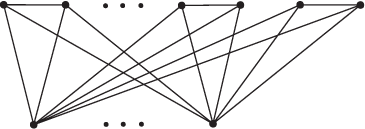}
  \caption{The graph $G=(k+1)K_{2}\vee kK_{1}$.}
\end{figure}

Obviously, the graph $G=(k+1)K_{2}\vee kK_{1}$  depicted in Fig. 4 is a connected graph with $\kappa(G)=k$,  $\alpha(G)=\kappa(G)+1$, and $|V(G)|=3\kappa(G)+2$. Since   $\omega(G-kK_{1})=k+1> k$, $G$ is not Hamiltonian. Suppose that $G$ has a 2-factor $C_{1}\cup C_{2}$ with two components. Denote  $C_{1}\cap kK_{1}=k_{1}K_{1}$, $C_{1}\cap (k+1)K_{2}=aK_{1}\cup bK_{2}$, $C_{2}\cap kK_{1}=k_{2}K_{1}$, $C_{2}\cap (k+1)K_{2}=aK_{1}\cup cK_{2}$, where $k_{1}+k_{2}=k$, $a+b+c=k+1$, $a,b,c$ are nonnegative integers, and $k_{1},k_{2}$ are positive integers. Since $C_{1}$ (resp. $C_{2}$) is a spanning cycle in $G[V(C_{1})]$  (resp. $G[V(C_{2})]$), we have $\omega(G[V(C_{1})]-k_{1}K_{1})=a+b\leq k_{1}$ and $\omega(G[V(C_{2})]-k_{2}K_{1})=a+c\leq k_{2}$. Then $2a+b+c=a+k+1\leq k_{1}+k_{2}=k$, a contradiction. Thus, $G$ has no 2-factor with at most two components and $G\notin \mathcal G$.  This means that the condition $|V(G)|\geq 3\kappa(G)+3$ in Theorem \ref{Th1}  is sharp.

Similarly, the graph $G=(k+2)K_{2}\vee kK_{1}$ is a connected graph with $\kappa(G)=k$, $\alpha(G)=\kappa(G)+2$, and $|V(G)|=3\kappa(G)+4$, while $G$ has no 2-factor with at most two components and $G\notin \mathcal G$. This means that  the condition $\alpha(G)\leq \kappa(G)+1$ in Theorem \ref{Th1} is best possible.

\section{Concluding remark}

In this paper, we main focus on the Chv\'{a}tal-Erd\H{o}s condition $\alpha(G)\leq \kappa(G)+1$, which   guaranteing that a connected graph has a 2-factor with at most two components. From the above discussion, we know that Theorem \ref{Th1},   Corollary  \ref{C2}, and Corollary \ref{C3} are all best possible. In addition, by Claim \ref{cl14}, we know that the result in Theorem \ref{Th6} still holds for a connected graph.

By Theorem \ref{Th4}, we know that every 2-connected graph $G$ with
$\alpha(G)= \kappa(G)+1$ can be covered by two cycles.  By Theorem \ref{Th3},  every  4-connected graph $G$  has a 2-factor with exactly two components if $\alpha(G)\leq \kappa(G)$ and $|V(G)|\geq 6$. By Corollary \ref{C2},  every connected graph $G$ of order at least $3\kappa(G) + 3$ and  minimum degree at least 2, with
$\alpha(G)\leq \kappa(G)+1$,  has a 2-factor with at most two components.
Motivated by these results, we pose two problems as follows.\\

\noindent\textbf{Problem $1$}: Characterize the conditions under which  a graph $G$ with $\alpha(G)= \kappa(G)+1$ and minimum degree at least 2  has a 2-factor with exactly two components.\\

\noindent\textbf{Problem $2$}: Does every Hamiltonian graph $G$ of order at least $3\kappa(G) + 3$ with $\alpha(G)\leq \kappa(G)+1$  have a 2-factor with exactly two components?\\


\noindent{\bf Declaration of competing interest}\\

There is no conflict of interest.

No experiments and no experimental data included in this article.\\

\noindent{\bf Data availability}\\

No data was used for the research described in the article.\\

\noindent{\bf Acknowledgements}\\

 We thank the National Institute of Education, Nanyang Technological University, where part of this research was performed. This work was supported by National Natural Science Foundation of
China (Nos. 12101126, 12131013 and 12571366), the Natural Science Foundation of Fujian Province (No. 2023J01539) and China Scholarship Council (No. 202409100010).

\end{CJK}

\end{document}